\documentclass[12pt,a4paper]{article}

\usepackage[margin=25mm]{geometry}
\usepackage[T1]{fontenc}
\usepackage[utf8]{inputenc}
\usepackage[english]{babel}
\usepackage{lmodern}
\usepackage{amsmath,amssymb,amsthm}
\usepackage{microtype}
\usepackage{cite}
\usepackage{needspace}
\usepackage[hidelinks]{hyperref}

\hypersetup{
  pdftitle={Two-point concentration for graph classes admitting sparse vertex deletion},
  pdfauthor={Margarita Akhmejanova and Vladislav Kozhevnikov},
  pdfsubject={Maximum induced subgraphs in dense binomial random graphs},
  pdfkeywords={random graphs, induced subgraphs, two-point concentration, vertex deletion},
  pdflang={en}
}

\numberwithin{equation}{section}
\theoremstyle{plain}
\newtheorem{theorem}{Theorem}
\newtheorem{lemma}[theorem]{Lemma}

\newtheorem{corollary}[theorem]{Corollary}
\newtheorem{example}[theorem]{Example}
\theoremstyle{remark}

\newcommand{\PP}{\mathcal P}
\newcommand{\Prob}{\mathbb P}
\newcommand{\E}{\mathbb E}
\newcommand{\Gnp}{G(n,p)}
\DeclareMathOperator{\Var}{Var}
\newcommand{\doi}[1]{\href{https://doi.org/#1}{\nolinkurl{doi:#1}}}

\title{Two-point concentration for graph classes\\admitting sparse vertex deletion}
\author{Margarita Akhmejanova\footnote{Matematiska institutionen,
Uppsala universitet, Box 480, 751 06, Uppsala, Sweden.\newline
Email: \href{mailto:mechmathrita@gmail.com}{\nolinkurl{mechmathrita@gmail.com}}.
Funded by the Wenner-Gren Foundation.}
\and Vladislav Kozhevnikov\footnote{Yandex.}}
\date{}

\begin{document}
\maketitle

\begin{abstract}
For fixed $p\in(0,1)$, we prove that the maximum number of vertices
in an induced subgraph of $G(n,p)$ from a given class is concentrated
on two consecutive integers, unifying and extending several earlier
results. We require that the class contain arbitrarily large graphs
and that every sufficiently large $k$-vertex member have a vertex
of degree $o(k)$ whose deletion leaves a graph in the same class.
This condition holds for empty graphs, trees, forests, planar
graphs, graphs of bounded maximum degree, degeneracy or treewidth,
graphs excluding a fixed complete bipartite subgraph (including
graphs containing no cycle of length four), and many other graph
classes.
\end{abstract}


\section{Introduction}\label{sec:introduction}

All graphs considered in the paper are simple, that is, undirected,
without loops or multiple edges. For a graph $H$, let $V(H)$ and $E(H)$
denote its vertex and edge sets, respectively. We call $|V(H)|$ the
size of $H$. For $U\subseteq V(G)$, $G[U]$ denotes the subgraph of $G$
induced by $U$. Graph classes are assumed to be invariant under
isomorphism.

Let $G(n,p)$ denote the binomial random graph on
$[n]:=\{1,\ldots,n\}$, in which each pair of vertices is joined by
an edge with probability $p$, independently of all other pairs.
We say that a property of $G(n,p)$ holds \emph{with high probability}
(\emph{whp}, for brevity) if its probability tends to $1$ as $n\to\infty$.
An integer-valued graph parameter is \emph{two-point concentrated} if whp
it belongs to a deterministic pair of consecutive integers. In this
paper, we study the maximum number of vertices in an induced subgraph
from a given graph class in $G(n,p)$. We prove that a simple condition
on vertex deletions implies
two-point concentration. We first give a brief overview of related
results, writing $q=(1-p)^{-1}$.

The first results of this type concern the independence number of
$G(n,p)$, or, equivalently, the clique number of $G(n,1-p)$.
Early estimates were obtained by Matula~\cite{Matula1970} and Grimmett
and McDiarmid~\cite{GrimmettMcDiarmid1975}. Frieze~\cite{Frieze1990}
obtained further asymptotic results in sparse random graphs. For
constant $p$, Matula~\cite{Matula1972,Matula1976} and, independently,
Bollob\'as and Erd\H{o}s~\cite{BollobasErdos1976} proved two-point
concentration. More precisely, for every constant $p\in(0,1)$ there is
a deterministic integer sequence $f(n)$ such that
\[
 f(n)=2\log_q n-2\log_q\log_q n+O(1)
\]
and whp
\[
 \alpha(G(n,p))\in\{f(n),f(n)+1\}.
\]
Thus, two-point concentration holds for the maximum size of an induced
edgeless graph.

Since then, two-point concentration has been proved for the maximum sizes
of induced subgraphs from several other classes. Fountoulakis, Kang
and McDiarmid~\cite{FKM2010} proved this for the $t$-stability number,
that is, the maximum size of an induced subgraph with maximum degree
at most $t$, when $t$ and $p$ are constant. In~\cite{FKM2014}, they
proved two-point concentration for induced subgraphs of bounded average
degree. For fixed $p\in(0,1)$, Kamaldinov, Skorkin and
Zhukovskii~\cite{KSZ2021} proved two-point concentration of the maximum
number of vertices $k$ in an induced subgraph of $G(n,p)$ with exactly
$e(k)$ edges, provided that $$e(k)<\varepsilon k^2$$ for sufficiently
small $\varepsilon>0$ and $t$ satisfies a suitable smoothness condition. In contrast, Balogh and Zhukovskii~\cite{BZ2022} proved that,
when $$e(k)=p\binom{k}{2}+O(k),$$ this maximum is not concentrated
on any fixed finite set of integers.

Dutta and Subramanian~\cite{DuttaSubramanian2018,DuttaSubramanian2023}
proved two-point concentration for the maximum sizes of induced paths
and cycles when $n^{-1/2}(\log n)^2\le p\le1-\varepsilon$,
where $\varepsilon>0$ is fixed. For induced matchings,
El Maftouhi and Marquez Gordones~\cite{EMG1994} proved concentration
in an interval of bounded length for constant $p$. Clark~\cite{Clark2001}
proved that, for constant $p$, the maximum number of edges in an induced matching is
concentrated on two consecutive integers. Cooley, Dragani\'c, Kang and
Sudakov~\cite{CDKS2021} obtained asymptotically optimal first-order
bounds for edge probabilities down to $p=C/n$, where $C$ is a
sufficiently large constant. Dragani\'c, Glock and
Krivelevich~\cite{DGK2022} extended these bounds to induced subgraphs
whose components are all isomorphic to an arbitrary fixed tree.

The maximum sizes of induced trees and forests are also two-point
concentrated when $p$ is constant. This was proved by Kamaldinov,
Skorkin and Zhukovskii~\cite{KSZ2021} for trees and by Krivoshapko
and Zhukovskii~\cite{KZ2021} for forests. In three
papers~\cite{AK2024,AK2024note,AKZ2026}, the present authors, jointly
with Zhukovskii in the third paper, obtained concentration in intervals
of length $o(1/p)$ for the maximum sizes of induced forests with and
without degree restrictions in sparse random graphs. We also proved two-point concentration for the maximum sizes of induced trees and
forests of bounded maximum degree for constant $p$.

Several recent results concern sparse random graphs. Bohman and
Hofstad~\cite{BH2024} proved two-point concentration of the independence
number for $p>n^{-2/3+\varepsilon}$, where $\varepsilon>0$ is constant.
In~\cite{BH2024nm}, they obtained the sharper sufficient condition
$p\gg(\log n/n)^{2/3}$ and proved two-point concentration in $G(n,m)$
for $m>n^{5/4+\varepsilon}$. For induced trees, Buitrago Oropeza~\cite{Oropeza2025} proved
two-point concentration when
\[
n^{-\frac{e-2}{3e-2}+\varepsilon}\le p,
\qquad p=o(1).
\]
Hofstad~\cite{Hofstad2026} extended this range to
\[
p\gg n^{-1/2}(\log n)^{3/2}.
\] Thus, for these
classes, two-point concentration also holds when $p\to0$, although the
known ranges of $p$ depend on the class.

There are also concentration results for random hypergraphs. Denisov
and Shabanov~\cite{DenisovShabanov2021} obtained concentration results
for their independence numbers. Vakhrushev~\cite{Vakhrushev2025} proved two-point concentration
of the weak independence number of random uniform hypergraphs for
some ranges of $p$. Demidovich and
Shabanov~\cite{DemidovichShabanov2022} obtained two-point concentration
results for the chromatic number of random hypergraphs.

We also mention results for dense induced subgraphs. Balister,
Bollob\'as, Sahasrabudhe and Veremyev~\cite{BBSV2019} proved two-point
concentration for the maximum size of an induced subgraph of density
at least $\gamma$, where $0<p<\gamma<1$ are constant. This class is
not covered by the condition in our main theorem. Bollob\'as and
Thomason~\cite{BT2000} developed a general first-order theory for the
maximum sizes of induced subgraphs with hereditary properties.
However, hereditary properties do not always give two-point
concentration. Indeed, Bohman, Michelen and Mubayi~\cite{BMM2026}
proved that the maximum size of an induced $K_r$-free subgraph of
$G(n,1/2)$ is two-point concentrated for $r=2,3$, but not for $r\ge4$.
Our result below gives a sufficient condition for two-point concentration
in terms of vertex deletions.

Finally, two-point concentration has been proved for other parameters
of random graphs. Wieland and Godbole~\cite{WG2001} proved it for
the domination number in dense random graphs. Glebov, Liebenau and
Szab\'o~\cite{GLS2015} extended the range to
$p\gg(\log n)^2/\sqrt n$, and Bohman, Warnke and
Zhu~\cite{BWZ2024} to $p\ge n^{-2/3+\varepsilon}$. Domination depends
on the edges between a chosen set and its complement, so it cannot
be described by a condition on the induced graph alone. For the
chromatic number, Shamir and Spencer~\cite{SS1987} developed
concentration results, and Alon and Krivelevich~\cite{AlonKrivelevich1997}
proved two-point concentration for $p\le n^{-1/2-\varepsilon}$.
Kargaltsev, Shabanov and Shaikheeva~\cite{KSS2019} identified two
explicit consecutive concentration values in certain ranges with
$p\le n^{-3/4-\delta}$ and $np\to\infty$, where $\delta>0$ is
constant. For dense random graphs, Heckel and Riordan~\cite{HR2023}
proved that the fluctuations of the chromatic number are much larger.

\subsection*{Main result}

For a graph class $\mathcal P$, let
\[
X_{\mathcal P}(G)=
\max
\{|S|:S\subseteq V(G),\ G[S]\in\mathcal P\}.
\]
If $S=\emptyset$, put $X_{\mathcal P}(G)=0.$
\begin{theorem}\label{thm:main}
Let $\mathcal P$ be a graph class containing graphs of arbitrarily
large size. Suppose that there is a function $d(k)=o(k)$ such that
every sufficiently large graph $H\in\mathcal P$ of size $k$ has a
vertex $v$ satisfying
\begin{gather}\label{eq:permitted}
H-v\in\mathcal P,
\qquad
\deg_H(v)\le d(k).
\end{gather}
Then, for every constant $p\in(0,1)$, there is a deterministic integer
sequence $f(n)$ such that
\[
\mathbb P\bigl(
X_{\mathcal P}(G(n,p))\in\{f(n),f(n)+1\}
\bigr)\longrightarrow1.
\]
Moreover, $f(n)\sim2\log_{1/(1-p)}n$.
\end{theorem}

Note that the assumption $d(k)=o(k)$ in Theorem~\ref{thm:main} cannot be relaxed to $d(k)=O(k)$
(see Example~\ref{ex:sharpness-d}) and that the conclusion can fail even when $p=p(n)\to0$ arbitrarily slowly
(see Example~\ref{ex:sharpness-p}).

A graph class is called \emph{hereditary} if it is closed under taking
induced subgraphs. For such classes, the condition on vertex deletions
follows from a bound on the number of edges.

\begin{corollary}\label{cor:hereditary}
Let $\mathcal P$ be a hereditary graph class containing graphs of
arbitrarily large size. Suppose that every $H\in\mathcal P$ of size
$k$ has $o(k^2)$ edges, uniformly over the class. Then the conclusion
of Theorem~\ref{thm:main} holds.
\end{corollary}

\begin{proof}
Every graph $H\in\mathcal P$ of size $k$ has a vertex of degree at
most $2|E(H)|/k=o(k)$. Since $\mathcal P$ is hereditary, deleting this
vertex leaves a graph in $\mathcal P$. Theorem~\ref{thm:main} applies.
\end{proof}

The condition of Theorem~\ref{thm:main} holds for empty graphs and
graphs of bounded maximum degree~\cite{Matula1972,BollobasErdos1976,FKM2010},
trees and forests with or without a fixed bound on the maximum
degree~\cite{KSZ2021,KZ2021,AKZ2026}, and paths and linear
forests~\cite{DuttaSubramanian2023}. It also holds for graphs of
bounded average degree and therefore gives the concentration result
of Fountoulakis, Kang and McDiarmid~\cite{FKM2014}.

\subsection{Idea of the proof and notation}

For the remainder of the paper, we fix $p\in(0,1)$.
Throughout the paper, we use natural logarithms unless a base is indicated.
All asymptotic relations are considered as $n\to\infty$, or as the specified
parameter tends to infinity. Constants implicit in asymptotic notation may
depend on $p$ and $\PP$. For integers $x\ge r\ge0$, let
$(x)_r=x(x-1)\cdots(x-r+1)$, where $(x)_0=1$. For a graph $F$ on a
given set of $j$ labelled vertices, let
\begin{equation}\label{eq:pi}
 \pi_j(F)=\Prob(G(j,p)=F)
 =p^{|E(F)|}(1-p)^{\binom j2-|E(F)|}
 =q^{-\binom j2}(pq)^{|E(F)|}.
\end{equation}
Note that $\pi_j(F)>0$ for every $F$.

The main step of the proof is to show that, for
$k=2\log_q n+o(\log n)$,
\begin{equation}\label{eq:strategy}
\Prob(X_{\PP}(\Gnp)\ge k)\ge\frac1{\log n}
\quad\Longrightarrow\quad
\Prob(X_{\PP}(\Gnp)\ge k-1)\longrightarrow1.
\end{equation}
We take an auxiliary copy of $\Gnp$ conditioned on the existence of an
induced subgraph in $\PP$ on $k$ vertices. We choose one such subgraph and
delete a vertex satisfying~\eqref{eq:permitted}. The remaining graph belongs
to $\PP$ and has $k-1$ vertices. Let $\nu_{k-1}$ be its distribution after
a uniform random labelling.

We then take an independent graph $G\sim\Gnp$. For each $(k-1)$-subset
$U\subseteq[n]$, let
\[
w_U=\frac{\nu_{k-1}(G[U])}{\pi_{k-1}(G[U])}.
\]
We apply the second moment method to
\[
T=\frac{1}{\binom{n}{k-1}}
\sum_{\substack{U\subseteq[n]\\|U|=k-1}}w_U.
\]
By the choice of the weights, $\E T=1$. Moreover, $T>0$ implies that $G$
contains an induced subgraph in $\PP$ on $k-1$ vertices. It therefore
suffices to prove that $\Var T=o(1)$.

To estimate the variance, we group pairs of $(k-1)$-subsets according to
the size of their intersection. For small and intermediate intersections,
we use the fact that every graph in $\PP$ on $k$ vertices has $o(k^2)$
edges. For intersections containing almost all $k-1$ vertices, we also
keep track of the deleted vertex. Its degree is $o(k)$, so there are only
$\exp(o(k))$ possible neighbourhoods. This gives the required bound for
the remaining terms in the second moment and proves~\eqref{eq:strategy}.

Finally, we choose the largest $k$ such that
$\Prob(X_{\PP}(\Gnp)\ge k)\ge1/\log n$. Together with the first-order
bounds,~\eqref{eq:strategy} gives two-point concentration.

The rest of the paper is organised as follows. In
Section~\ref{sec:preliminaries}, we prove several elementary estimates.
In Section~\ref{sec:location}, we determine the first-order asymptotics of
$X_{\PP}(\Gnp)$. Section~\ref{sec:amplification} is devoted to the proof of
the main probabilistic lemma. Finally, in Section~\ref{sec:quantile}, we
complete the proof of Theorem~\ref{thm:main}.

\section{Preliminaries}\label{sec:preliminaries}

Under the assumptions of Theorem~\ref{thm:main}, fix an integer
$k_0\ge2$ such that every graph $H\in\PP$ on $k\ge k_0$ vertices has a
vertex $v\in V(H)$ satisfying~\eqref{eq:permitted}. Such an integer exists
by the assumptions of the theorem.

Since $\PP$ contains graphs of arbitrarily large size, repeated deletions
of vertices satisfying~\eqref{eq:permitted} give a member of $\PP$ on
every number $k\ge k_0-1$ of vertices. The same argument shows that, for
every graph $G$ and every such $k$,
\begin{equation}\label{eq:exactorder}
 X_{\PP}(G)\ge k
 \quad\Longleftrightarrow\quad
 \exists U\subseteq V(G):\ |U|=k,\ G[U]\in\PP.
\end{equation}

\begin{lemma}\label{lem:chains}
Let $H\in\PP$ have $k\ge k_0-1$ vertices. Then
\begin{equation}\label{eq:sparsebound}
 |E(H)|\le s(k):=\binom{k_0-1}{2}+\sum_{t=k_0}^{k}d(t)=o(k^2).
\end{equation}
The sum is empty when $k=k_0-1$.
\end{lemma}

\begin{proof}
We successively delete vertices satisfying~\eqref{eq:permitted} until
$k_0-1$ vertices remain. At the step from $t$ to $t-1$ vertices, we remove
at most $d(t)$ edges. The remaining graph has at most
$\binom{k_0-1}{2}$ edges, which proves the inequality. Finally,
$\sum_{t=k_0}^{k}d(t)=o(k^2)$ since $d(t)=o(t)$.
\end{proof}

The function $s$ is non-decreasing, and the bound holds uniformly over
all graphs in $\PP$ with the given number of vertices.

\subsection{Entropy estimates}

Let
\[
 \mathcal H(x)=-x\log x-(1-x)\log(1-x),\qquad 0<x<1,
\]
and set $\mathcal H(0)=\mathcal H(1)=0$. We will use the following
standard estimates.

\begin{lemma}\label{lem:entropy}
Let $M\ge1$ and $0\le r\le M$ be integers. Then
\begin{equation}\label{eq:binentropy}
 \binom Mr\le\exp\bigl(M\mathcal H(r/M)\bigr).
\end{equation}
If $r\le M/2$ and $z>0$, then
\begin{equation}\label{eq:sumentropy}
 \sum_{h=0}^{r}\binom Mh z^h
 \le(r+1)\exp\bigl(M\mathcal H(r/M)\bigr)\max\{1,z\}^{r}.
\end{equation}
In particular, for fixed $z$ and $r=o(M)$,
\begin{equation}\label{eq:entropyasymp}
 \log\binom Mr=o(M),\qquad
 \log\left(\sum_{h=0}^{r}\binom Mh z^h\right)=o(M).
\end{equation}
\end{lemma}

\begin{proof}
If $0<r<M$, let $x=r/M$. Since
$\binom Mr x^r(1-x)^{M-r}\le(x+1-x)^M=1$, we get~\eqref{eq:binentropy}.
For $r=0$ and $r=M$, the inequality is immediate.
The function $\mathcal H$ is increasing on $[0,1/2]$ since
$\mathcal H'(x)=\log((1-x)/x)>0$ for $0<x<1/2$. Thus every summand
in~\eqref{eq:sumentropy} is at most
$\exp(M\mathcal H(r/M))\max\{1,z\}^{r}$, which proves the second
inequality. Finally, if $r=o(M)$, then $\mathcal H(r/M)=o(1)$,
$\log(r+1)=o(M)$ and $r\log\max\{1,z\}=o(M)$, giving~\eqref{eq:entropyasymp}.
\end{proof}

\begin{lemma}\label{lem:overlap}
Let $U,V$ be chosen independently and uniformly from the $m$-subsets of
$[n]$, where $n>m$, and let
\begin{equation}\label{eq:Aj}
 A_j^{(m)}=\Prob(|U\cap V|=j)
 =\frac{\binom mj\binom{n-m}{m-j}}{\binom nm},\qquad 0\le j\le m.
\end{equation}
Then
\begin{equation}\label{eq:overlapupper}
 A_j^{(m)}\le\left(\frac{m^2}{n-m}\right)^j.
\end{equation}
Moreover, if $n\ge2m$, then
\begin{equation}\label{eq:overlapidentity}
 A_j^{(m)}\binom nj
 =\binom mj^2\frac{\binom{n-m}{m-j}}{\binom{n-j}{m-j}}
 \le\binom mj^2.
\end{equation}
\end{lemma}

\begin{proof}
The formula for $A_j^{(m)}$ follows by first choosing $U\cap V$ and then
choosing the vertices of $V$ outside $U$. Let $J=|U\cap V|$. For $j\ge1$,
\[
 A_j^{(m)}\le\Prob(J\ge j)\le\E\binom Jj
 =\binom mj\frac{(m)_j}{(n)_j}
 \le\left(\frac{m^2}{n-m}\right)^j.
\]
For $j=0$, the bound is $A_0^{(m)}\le1$.
To prove~\eqref{eq:overlapidentity}, we use the identity
$\binom nm\binom mj=\binom nj\binom{n-j}{m-j}$.
The ratio in~\eqref{eq:overlapidentity} is at most one since
$n-m\le n-j$. Its denominator is positive when $n\ge2m$.
\end{proof}

\section{First-order bounds}\label{sec:location}

To obtain the lower bound for $X_{\PP}(\Gnp)$, it is enough to find
an induced copy of one graph from $\PP$ of the required size.
We start with a lemma on induced copies of a given sparse graph.

\begin{lemma}\label{lem:embedding}
Fix $\eta\in(0,2)$ and set $k=\lfloor(2-\eta)\log_q n\rfloor$.
Let $H_n$ be a deterministic graph on $k$ vertices with
$|E(H_n)|=o(k^2)$. Then $\Gnp$ contains an induced copy of $H_n$ whp.
\end{lemma}

\begin{proof}
Let $H=H_n$ be labelled by $[k]$, and let $Y$ denote the number of
injections $\phi:[k]\to[n]$ which induce copies of $H$.
Then $\E Y=(n)_k\pi_k(H)>0$.
Consider two injections whose images intersect in $j$ vertices.
If they require different graphs on their intersection, the probability
that both induce copies of $H$ is zero. Otherwise, let $t$ be the number
of edges in the common graph. The ratio of the joint probability to
the product of the two probabilities is
\[
 p^{-t}(1-p)^{-\binom j2+t}
 =q^{\binom j2}(pq)^{-t}
 \le\exp\left(\binom j2\log q+
 \min\left\{|E(H)|,\binom j2\right\}\log(1/p)\right).
\]
For $j=0,1$, the corresponding indicators are independent, since they
depend on disjoint sets of edges. The image of a uniform injection is
a uniform $k$-subset of $[n]$. Therefore,
\begin{equation}\label{eq:embeddingvariance}
 \frac{\Var Y}{(\E Y)^2}
 \le\sum_{j=2}^{k} A_j^{(k)}
 \exp\left(\binom j2\log q+
 \min\left\{|E(H)|,\binom j2\right\}\log(1/p)\right).
\end{equation}
By~\eqref{eq:overlapupper}, the logarithm of each summand on the
right-hand side is at most
\begin{equation}\label{eq:embeddingexponent}
 -j\log n+2j\log k+o(j)+\frac{j^2}{2}\log q
 +\min\left\{|E(H)|,\binom j2\right\}\log(1/p).
\end{equation}
The error term obtained from $\log(n-k)=\log n+o(1)$ is uniform in $j$.

We estimate the sum in~\eqref{eq:embeddingvariance} in two ranges.
For $2\le j\le\log n/(2\log(q/p))$, the quadratic terms
in~\eqref{eq:embeddingexponent} satisfy
\[
 \frac{j^2}{2}\log(q/p)\le\frac j4\log n.
\]
Since $\log k=o(\log n)$, the exponent is at most $-j\log n/2$
for all sufficiently large $n$, uniformly over this range. Thus, the
sum over these $j$ is at most
\[
 \sum_{j=2}^{\infty}n^{-j/2}
 =\frac{n^{-1}}{1-n^{-1/2}}=o(1).
\]
Now let $j>\log n/(2\log(q/p))$. Since $j\le k$, we have
$j^2\log q/2\le(1-\eta/2)j\log n$. Also,
\[
 2j\log k+o(j)+|E(H)|\log(1/p)=o(j\log n)
\]
uniformly, since $j$ is bounded below by a positive constant times
$\log n$ and $|E(H)|=o((\log n)^2)$. Thus the exponent is at most
$-\eta j\log n/4$. There are $O(\log n)$ such terms, each at most
$\exp(-C_{\eta,p}(\log n)^2)$ for a constant $C_{\eta,p}>0$.
It follows that the right-hand side of~\eqref{eq:embeddingvariance}
is $o(1)$. Applying Chebyshev's inequality, we get
$\Prob(Y=0)\le\Var Y/(\E Y)^2=o(1)$.
\end{proof}

\subsection{First and second moment bounds}

\begin{lemma}\label{lem:coarse}
Let the assumptions of Theorem~\ref{thm:main} hold. Fix $\eta\in(0,2)$
and set
\[
 k_- =\lfloor(2-\eta)\log_q n\rfloor,
 \qquad k_+=\lceil(2+\eta)\log_q n\rceil.
\]
Then
\begin{equation}\label{eq:coarsebounds}
 \Prob(X_{\PP}(\Gnp)\ge k_-)\longrightarrow1,
 \qquad
 \Prob(X_{\PP}(\Gnp)\ge k_+)
 \le\exp\bigl(-\gamma_\eta(\log n)^2\bigr)
\end{equation}
for all sufficiently large $n$, where $\gamma_\eta>0$.
In particular,
\begin{equation}\label{eq:firstorder}
 \frac{X_{\PP}(\Gnp)}{2\log_q n}\xrightarrow{\Prob}1.
\end{equation}
\end{lemma}

\begin{proof}
By the observation preceding Lemma~\ref{lem:chains}, for all large
$n$ there exists a deterministic graph $H_n\in\PP$ of size $k_-$.
Lemma~\ref{lem:chains} gives $|E(H_n)|\le s(k_-)=o(k_-^2)$.
The lower bound now follows from Lemma~\ref{lem:embedding}.

To prove the upper bound, let $\PP_k$ be the set of graphs from $\PP$
on $[k]$, and let $Z_k$ denote the number of $k$-subsets $U\subseteq[n]$
such that $G[U]\in\PP$. Then
\begin{equation}\label{eq:EZ}
 \E Z_k=\binom nk q^{-\binom k2}
 \sum_{H\in\PP_k}(pq)^{|E(H)|}.
\end{equation}
Since $s(k)=o(k^2)$, for sufficiently large $k$ we have
$\lfloor s(k)\rfloor\le\binom k2/2$. Lemma~\ref{lem:entropy} therefore gives
\[
 \sum_{H\in\PP_k}(pq)^{|E(H)|}
 \le\sum_{h=0}^{\lfloor s(k)\rfloor}
 \binom{\binom k2}{h}(pq)^h
 =\exp(o(k^2)).
\]
Since $\binom nk\le n^k$, we have
\begin{equation}\label{eq:firstmomentlog}
 \log\E Z_k
 \le k\log n-\frac{k(k-1)}2\log q+o(k^2)
 =k\log n-\frac{k^2}2\log q+o(k^2).
\end{equation}
At $k=k_+$, the right-hand side equals
\[
 -\frac{\eta(2+\eta)}{2\log q}(\log n)^2+o((\log n)^2).
\]
Therefore, we may take $\gamma_\eta=\eta(2+\eta)/(4\log q)$
for all sufficiently large $n$.
By~\eqref{eq:exactorder} and Markov's inequality,
\[
 \Prob(X_{\PP}\ge k_+)=\Prob(Z_{k_+}>0)
 \le\E Z_{k_+}\le\exp\bigl(-\gamma_\eta(\log n)^2\bigr).
\]
Since $\eta>0$ is arbitrary, the two bounds imply~\eqref{eq:firstorder}.
\end{proof}

We will use the upper bound in~\eqref{eq:coarsebounds} in
Section~\ref{sec:quantile} to choose a quantile with tail probability
$1/\log n$.

\section{Concentration via vertex deletion}\label{sec:amplification}

In this section, we prove the following implication, which is the main
step in the proof of Theorem~\ref{thm:main}.

\begin{lemma}\label{lem:amplification}
Let $k=k(n)$ be a deterministic integer sequence such that
\begin{equation}\label{eq:knear}
 k=2\log_q n+o(\log n).
\end{equation}
If
\begin{equation}\label{eq:epscondition}
 \varepsilon:=\Prob(X_{\PP}(\Gnp)\ge k)\ge\frac1{\log n}
\end{equation}
for all sufficiently large $n$, then
\begin{equation}\label{eq:amplificationconclusion}
 \Prob(X_{\PP}(\Gnp)\ge k-1)\longrightarrow1.
\end{equation}
\end{lemma}

Let $n$ be large enough that $k\ge k_0$, and put $m=k-1$ and
$D=\lfloor d(k)\rfloor$. Then
\begin{equation}\label{eq:mD}
 m=2\log_q n+o(\log n),\qquad D=o(m),\qquad s(k)=o(m^2).
\end{equation}
All distributions and weights considered below may depend on $n$.

\subsection{A weighted count of induced subgraphs}

Let $A_k=\{X_{\PP}(\Gnp)\ge k\}$. We sample $\widetilde G$ from
$\Gnp$ conditioned on $A_k$. By~\eqref{eq:exactorder}, there is a set
$S$ of exactly $k$ vertices such that $\widetilde G[S]\in\PP$.
We choose $S$ by an arbitrary fixed rule and then use another fixed
rule to choose a vertex $v\in S$ satisfying~\eqref{eq:permitted}.

Independently of these choices, we label the vertices of
$S\setminus\{v\}$ by a uniformly chosen bijection from $[m]$.
Let $\nu_m$ be the distribution of the resulting labelled graph after
deletion, and let $\nu_j$ be its restriction to the first $j$ labels,
where $0\le j\le m$. The random labelling implies that each $\nu_j$
is invariant under permutations of its labels. In particular, any fixed
set of $j$ labels gives the same distribution $\nu_j$ after relabelling
by $[j]$. The distribution $\nu_m$ is supported on $\PP_m$.
Graphs with fewer retained vertices need not belong to $\PP$, but
Lemma~\ref{lem:chains} gives
\begin{equation}\label{eq:restrictionedges}
 \nu_j(F)>0\quad\Longrightarrow\quad
 |E(F)|\le\min\left\{s(k),\binom j2\right\}.
\end{equation}
The probabilities $\nu_j(F)$ are deterministic.

We now take an independent graph $G\sim\Gnp$ without conditioning.
For each $m$-subset $U\subseteq[n]$, define
\begin{equation}\label{eq:weights}
 w_U(G)=\frac{\nu_m(G[U])}{\pi_m(G[U])},
 \qquad T(G)=\frac1{\binom nm}\sum_{|U|=m}w_U(G).
\end{equation}
By permutation invariance, the ratio does not depend on the bijection
used to identify $U$ with $[m]$. Every weight is nonnegative, and
\[
 \E w_U=\sum_{H\text{ on }[m]}\pi_m(H)
 \frac{\nu_m(H)}{\pi_m(H)}=1.
\]
Furthermore, if $w_U(G)>0$, then $G[U]\in\PP$, since $\nu_m$ is
supported on $\PP_m$. Consequently,
\begin{equation}\label{eq:ET}
 \E T=1,\qquad
 \{T>0\}\subseteq\{X_{\PP}(G)\ge m\}.
\end{equation}
The expectation in~\eqref{eq:ET} and all further moments of $T$ are
taken with respect to this unconditioned graph $G$.

\subsection{The second moment}

We estimate $\E T^2$ by grouping pairs of $m$-subsets according to
the number of common vertices. Fix $U,V$ with $I=U\cap V$ having
$j$ vertices, and label both sets so that the vertices of $I$ receive
the same first $j$ labels. Given $G[I]=F$, a labelled graph $H$ on
$U$ with $H[[j]]=F$ has probability $\pi_m(H)/\pi_j(F)$. Hence,
\begin{align*}
 \E[w_U\mid G[I]=F]
 &=\sum_{H:\ H[[j]]=F}\frac{\pi_m(H)}{\pi_j(F)}
                    \frac{\nu_m(H)}{\pi_m(H)}\\
 &=\frac{\nu_j(F)}{\pi_j(F)}.
\end{align*}
The same equality holds for $w_V$. Once $G[I]$ is fixed, the remaining
edge variables determining $w_U$ and $w_V$ are disjoint and independent.
Therefore,
\begin{equation}\label{eq:Rj}
 \E[w_Uw_V]=R_j,
 \qquad R_j:=\sum_{F\text{ on }[j]}\frac{\nu_j(F)^2}{\pi_j(F)}.
\end{equation}
Summing over all pairs of $m$-subsets, we get
\begin{equation}\label{eq:ET2}
 \E T^2=\sum_{j=0}^{m}A_j^{(m)}R_j.
\end{equation}
For $j=0,1$, there is only one graph on $[j]$, so $R_0=R_1=1$.
Since $\sum_j A_j^{(m)}=1$, it follows that
\begin{equation}\label{eq:variancebound}
 \Var T=\sum_{j=2}^{m}A_j^{(m)}(R_j-1)
 \le\sum_{j=2}^{m}A_j^{(m)}R_j.
\end{equation}

We use two bounds for $R_j$. First, by~\eqref{eq:restrictionedges},
every graph with positive $\nu_j$-probability has at most
$\min\{s(k),\binom j2\}$ edges. Hence,
\begin{equation}\label{eq:sparseR}
 R_j\le\exp\left(\binom j2\log q+
 \min\left\{s(k),\binom j2\right\}\log(1/p)\right),
\end{equation}
since $\sum_F\nu_j(F)^2\le1$ and
$\pi_j(F)^{-1}=q^{\binom j2}(pq)^{-|E(F)|}\le q^{\binom j2}p^{-|E(F)|}$.
For $j$ close to $m$, we also use the degree bound on the deleted
vertex, as follows.

\begin{lemma}\label{lem:root}
For $0\le j\le m$, set
\begin{equation}\label{eq:Bj}
 B_j=\sum_{h=0}^{\min\{D,j\}}\binom jh(pq)^h.
\end{equation}
Then, for every graph $F$ on $[j]$,
\begin{equation}\label{eq:marginalroot}
 \frac{\nu_j(F)}{\pi_j(F)}
 \le\frac{(n-j)\binom nj}{\varepsilon}\,q^{-j}B_j.
\end{equation}
In particular,
\begin{equation}\label{eq:rootR}
 R_j\le\frac{(n-j)\binom nj}{\varepsilon}\,q^{-j}B_j.
\end{equation}
Moreover, $B_j\le B_m$ and $\log B_m=o(m)$.
\end{lemma}

\begin{proof}
We consider again $\widetilde G$ and the selected vertex $v$, which
we label by $*$. Let $\rho_j(J)$ be the probability that $v$ together
with the first $j$ labels after deletion induces the labelled graph
$J$ on $\{*,1,\ldots,j\}$.

In an unconditioned graph $\Gnp$, let $N_J$ count ordered
$(j+1)$-tuples of distinct vertices inducing $J$, with the first entry
corresponding to $*$. Then
\[
 \E N_J=(n)_{j+1}\pi_{j+1}(J).
\]
Given $\widetilde G$, $S$ and $v$, the first $j$ labels after deletion
form a uniform ordered $j$-tuple from $S\setminus\{v\}$.
The number of such tuples inducing $J$ together with $v$ is at most
the total number of ordered tuples inducing $J$ in $\widetilde G$.
Since $\widetilde G$ has the distribution of $\Gnp$ conditioned on
$A_k$, taking expectations gives
\begin{equation}\label{eq:rootlaw}
\begin{aligned}
 \rho_j(J)
 &\le\frac{\E[N_J\mid A_k]}{(m)_j}
 \le\frac{\E N_J}{\varepsilon(m)_j}\\
 &\le\frac{(n)_{j+1}}{\varepsilon j!}\pi_{j+1}(J)
 =\frac{(n-j)\binom nj}{\varepsilon}\pi_{j+1}(J).
\end{aligned}
\end{equation}
Here $(m)_j\ge j!$, since $j\le m$. The argument uses only the
uniform random labelling; the rules for choosing $S$ and $v$ are
arbitrary and fixed in advance.

Fix a graph $F$ on $[j]$. For $B\subseteq[j]$, let $J(F,B)$ be the
graph obtained by adding to $F$ the vertex $*$ with neighbourhood
$B$. Then
\begin{equation}\label{eq:rootpi}
 \pi_{j+1}(J(F,B))
 =\pi_j(F)p^{|B|}(1-p)^{j-|B|}
 =\pi_j(F)q^{-j}(pq)^{|B|}.
\end{equation}
The degree of $v$ in the selected graph $\widetilde G[S]$ is at most
$D$. Thus $\rho_j(J(F,B))=0$ for $|B|>D$.
Summing~\eqref{eq:rootlaw} over all neighbourhoods of size at most
$D$ and applying~\eqref{eq:rootpi}, we obtain~\eqref{eq:marginalroot}.
Since
\[
 R_j=\sum_F\nu_j(F)\frac{\nu_j(F)}{\pi_j(F)},
 \qquad\sum_F\nu_j(F)=1,
\]
the same bound holds for $R_j$.

For $j\le m$, every term in $B_j$ is bounded by the corresponding
term in $B_m$. Also, $D=o(m)$, so $D\le m/2$ for sufficiently
large $n$. By Lemma~\ref{lem:entropy},
\begin{equation}\label{eq:Bm}
 \log B_m\le\log(D+1)+m\mathcal H(D/m)
             +D\log\max\{1,pq\}=o(m).
\end{equation}
This also holds for $D=0$, when $B_m=1$.
\end{proof}

\subsection{The three overlap ranges}

It remains to show that the sum on the right-hand side
of~\eqref{eq:variancebound} is $o(1)$. We choose a positive integer
$r_0$ such that
\begin{equation}\label{eq:cutoffproperties}
 r_0=o(m),\qquad
 |m-2\log_q n|+\log m+\frac{s(k)}m+1=o(r_0).
\end{equation}
Such a choice is possible by~\eqref{eq:mD}, since the sum on the left
is $o(m)$. Set
\[
 \theta=\frac{\log q}{8\log(q/p)}\in(0,1/8).
\]
For sufficiently large $n$, we have $m-r_0>\theta m$.
We consider the following three cases. All estimates below are uniform
over the indicated values of $j$.

\paragraph{Small overlaps: $2\le j\le\theta m$.}
By~\eqref{eq:overlapupper} and~\eqref{eq:sparseR},
\begin{align*}
 \log(A_j^{(m)}R_j)
 &\le j\bigl(2\log m-\log(n-m)\bigr)
       +\binom j2\log(q/p)\\
 &\le j\left(-\log n+2\log m+o(1)+\frac j2\log(q/p)\right).
\end{align*}
Since $m/\log n\to2/\log q$ and
$\theta\log(q/p)/\log q=1/8$, the expression in parentheses is at
most $-\log n/2$ for sufficiently large $n$. Hence,
\begin{equation}\label{eq:rangeI}
 \sum_{2\le j\le\theta m}A_j^{(m)}R_j
 \le\sum_{j=2}^{\infty}n^{-j/2}=o(1).
\end{equation}

\paragraph{Intermediate overlaps: $\theta m<j\le m-r_0$.}
Using~\eqref{eq:sparseR} and $j>\theta m$, we get
\begin{align*}
 \log(A_j^{(m)}R_j)
 &\le j(-\log n+2\log m+o(1))
       +\frac{j^2}{2}\log q+s(k)\log(1/p)\\
 &\le j\left[-\log n+2\log m+o(1)
       +\frac{m-r_0}{2}\log q
       +\frac{s(k)\log(1/p)}{\theta m}\right]\\
 &=j\left[-\frac{r_0}{2}\log q
       +\frac{m-2\log_q n}{2}\log q+2\log m
       +\frac{s(k)\log(1/p)}{\theta m}+o(1)\right].
\end{align*}
By the choice of $r_0$ in~\eqref{eq:cutoffproperties}, all terms in
brackets except $-r_0\log q/2$ are $o(r_0)$.
Thus the expression in brackets is at most $-r_0\log q/4$ for
sufficiently large $n$, and
\begin{equation}\label{eq:rangeII}
 \sum_{\theta m<j\le m-r_0}A_j^{(m)}R_j
 \le m\exp(-\theta mr_0\log q/4)=o(1).
\end{equation}

\paragraph{Nearly complete overlaps: $m-r_0<j\le m$.}
Write $j=m-r$, where $0\le r<r_0$. Applying~\eqref{eq:rootR}
and~\eqref{eq:overlapidentity}, we obtain
\begin{align}
 A_j^{(m)}R_j
 &\le\frac n\varepsilon q^{-j}B_j A_j^{(m)}\binom nj\notag\\
 &\le\frac n\varepsilon q^{-j}B_m\binom mr^2.
 \label{eq:nearbound}
\end{align}
Here~\eqref{eq:overlapidentity} cancels the factor $\binom nj$.
Since $j=m-r=2\log_q n+o(\log n)$ uniformly in this case,
\[
 \log(nq^{-j})=\log n-j\log q=-\log n+o(\log n).
\]
Also, $\log B_m=o(m)=o(\log n)$ by~\eqref{eq:Bm}.
Since $r_0=o(m)$, we have $r_0\le m/2$ for sufficiently large $n$,
and~\eqref{eq:binentropy} gives
\[
 \log\binom mr\le m\mathcal H(r/m)
 \le m\mathcal H(r_0/m)=o(m)=o(\log n).
\]
Finally, $\log(\varepsilon^{-1})\le\log\log n=o(\log n)$.
The right-hand side of~\eqref{eq:nearbound} is therefore
$n^{-1+o(1)}$, uniformly in this case. There are at most
$m=O(\log n)$ terms, so
\begin{equation}\label{eq:rangeIII}
 \sum_{m-r_0<j\le m}A_j^{(m)}R_j=o(1).
\end{equation}

Combining~\eqref{eq:rangeI}, \eqref{eq:rangeII}
and~\eqref{eq:rangeIII}, we get $\Var T=o(1)$. Since $\E T=1$,
Chebyshev's inequality and~\eqref{eq:ET} imply
\[
 \Prob(X_{\PP}(G)<k-1)\le\Prob(T=0)\le\Var T=o(1),
\]
which proves Lemma~\ref{lem:amplification}.
The argument uses only $d(k)=o(k)$ and $s(k)=o(k^2)$, with no
prescribed rate of convergence. The entropy estimates
in~\eqref{eq:Bm} and~\eqref{eq:nearbound} are needed here, since
bounds of the form $O(D\log m)$ and $O(r\log m)$ need not be $o(m)$.

\section{Proof of the main theorem}\label{sec:quantile}

\begin{proof}[Proof of Theorem~\ref{thm:main}]
Let $n$ be sufficiently large that $\log n>1$, and let
\begin{equation}\label{eq:quantile}
 k_n=\max\left\{t\in\{0,1,\ldots,n\}:
       \Prob(X_{\PP}(\Gnp)\ge t)\ge\frac1{\log n}\right\}.
\end{equation}
The set in~\eqref{eq:quantile} is nonempty since $X_{\PP}\ge0$.

Fix $\eta\in(0,2)$. By Lemma~\ref{lem:coarse}, the probability
$\Prob(X_{\PP}(\Gnp)\ge k_-)$ tends to $1$, where
$k_- =\lfloor(2-\eta)\log_q n\rfloor$. For
$k_+=\lceil(2+\eta)\log_q n\rceil$, the same lemma gives
\[
 \Prob(X_{\PP}(\Gnp)\ge k_+)
 \le\exp\bigl(-\gamma_\eta(\log n)^2\bigr)<\frac1{\log n}
\]
for all sufficiently large $n$. Therefore, $k_-\le k_n<k_+$.
Since $\eta>0$ is arbitrary, we get
\begin{equation}\label{eq:quantilelocation}
 k_n=2\log_q n+o(\log n).
\end{equation}
By the definition of $k_n$, $\Prob(X_{\PP}\ge k_n)\ge1/\log n$.
Applying Lemma~\ref{lem:amplification}, we obtain
\[
 \Prob(X_{\PP}<k_n-1)=o(1).
\]
Also, $k_n<n$ for all sufficiently large $n$ by~\eqref{eq:quantilelocation}.
The choice of $k_n$ in~\eqref{eq:quantile} then gives
$\Prob(X_{\PP}\ge k_n+1)<1/\log n$. Since $X_{\PP}$ is integer-valued,
\[
 \Prob\bigl(X_{\PP}\notin\{k_n-1,k_n\}\bigr)
 \le o(1)+\frac1{\log n}\longrightarrow0.
\]
Let $f(n)=k_n-1$ for all sufficiently large $n$, and define $f(n)$
arbitrarily for the remaining finitely many $n$. Then whp
$X_{\PP}(\Gnp)\in\{f(n),f(n)+1\}$, and~\eqref{eq:quantilelocation}
gives $f(n)\sim2\log_q n$, completing the proof.
\end{proof}

\subsection{Graphs of bounded average degree}

We now check that Theorem~\ref{thm:main} applies to graphs of bounded
average degree. The average degree may increase when a vertex is deleted,
so we need to choose the vertex appropriately. The following lemma gives
a suitable vertex of degree $O(\sqrt{k})$ in every sufficiently large
graph on $k$ vertices.

\begin{lemma}\label{lem:average-deletion}
For every constant $t\ge0$, there are constants $K_t\ge2$ and
$M_t<\infty$ such that every graph $H$ on $k\ge K_t$ vertices with
$|E(H)|\le tk/2$ has a vertex $v$ satisfying
\begin{equation}\label{eq:averagedeletion}
 \deg_H(v)\le M_t\sqrt{k},\qquad |E(H-v)|\le\frac{t(k-1)}2.
\end{equation}
For $t=0$, the degree bound may be taken to be zero.
\end{lemma}

\begin{proof}
The case $t=0$ is immediate. Let $t>0$, write $e=|E(H)|$, and let
\[
 \tau=e-\frac{t(k-1)}2.
\]
A vertex $v$ satisfies $|E(H-v)|\le t(k-1)/2$ if and only if
$\deg_H(v)\ge\tau$.
If $\tau\le0$, every vertex satisfies this condition. We may then choose
a vertex of minimum degree, which has degree at most $2e/k\le t$.

Now suppose that $\tau>0$. Since $e\le tk/2$, we have $\tau\le t/2$.
Let
\[
 S=\{v\in V(H):\deg_H(v)\ge\tau\},\quad h=|S|,
 \quad \ell=\left\lceil\frac t2\right\rceil-1,
 \quad \gamma=\frac t2-\ell>0.
\]
Every vertex outside $S$ has integer degree less than $\tau\le t/2$,
and hence at most $\ell$. Counting the edges inside $S$ and those with
at least one endpoint outside $S$, we get
\[
 e\le |E(H[S])|+\sum_{u\notin S}\deg_H(u)
 \le\binom h2+\ell(k-h)\le\frac{h^2}{2}+\ell k.
\]
Since $e>t(k-1)/2$, it follows that
\[
 \frac{h^2}{2}>\gamma k-\frac t2.
\]
For $k\ge t/\gamma$, this implies $h>\sqrt{\gamma k}$.
The sum of all degrees is at most $tk$, so there is a vertex $v\in S$
such that
\[
 \deg_H(v)\le\frac{\sum_{u\in S}\deg_H(u)}h
 \le\frac{tk}{h}\le\frac t{\sqrt\gamma}\sqrt{k}.
\]
Since $v\in S$, deleting $v$ leaves at most $t(k-1)/2$ edges.
Thus, we may take
\[
 K_t=\left\lceil\max\{2,t/\gamma\}\right\rceil,
 \qquad M_t=t\max\{1,\gamma^{-1/2}\}
\]
in both cases. For $t=0$, take $K_0=2$ and $M_0=0$.
\end{proof}

\begin{corollary}\label{cor:average-degree}
For every constant $t\ge0$, let $\PP_t=\{H:|E(H)|\le t|V(H)|/2\}$.
Then $X_{\PP_t}(\Gnp)$ is concentrated on two consecutive integers
with leading term $2\log_q n$.
\end{corollary}

\begin{proof}
The class $\PP_t$ is invariant under isomorphism and contains all
empty graphs. By Lemma~\ref{lem:average-deletion}, it satisfies
the condition of Theorem~\ref{thm:main} with
$d(k)=M_t\sqrt{k}=o(k)$ for all sufficiently large $k$.
The result follows from that theorem.
\end{proof}

This gives the qualitative concentration result of Fountoulakis, Kang
and McDiarmid~\cite{FKM2014} for constant $t$ and $p$.
Note that $\PP_t$ is not hereditary when $t>0$. Indeed, a sufficiently
large clique together with sufficiently many isolated vertices belongs
to $\PP_t$, while the induced clique has average degree greater than $t$.
A bound on the maximum average degree requires the same inequality
in every nonempty subgraph. This defines a hereditary class, to which
Corollary~\ref{cor:hereditary} applies directly.

\begin{example}[Failure for arbitrarily slowly vanishing $p$]
\label{ex:sharpness-p}
For every positive function $a(n)\to0$, there exist integers
$0=b_0<b_1<\cdots$ and a non-increasing function $p(n)\to0$,
with $p(n)\ge a(n)$ eventually, such that the fixed class
\[
\mathcal B=
\bigcup_{r=1}^{\infty}
\{K_r\sqcup(k-r)K_1:\ b_{r-1}<k\le b_r\}
\]
satisfies the deletion condition with $d(k)=o(k)$, but
$X_{\mathcal B}(G(n,p(n)))$ is not two-point concentrated.
\end{example}

\begin{proof}
The orders $b_{r-1}+1,\ldots,b_r$ correspond to clique
order $r$, and larger orders require larger cliques.
Thus, if $G$ contains no $K_r$, then
$X_{\mathcal B}(G)\le b_{r-1}$.
Conversely, an $r$-clique together with $b_r-r$ mutually
nonadjacent vertices having no neighbours in the clique
induces a graph in $\mathcal B$ of order $b_r$.
We arrange that these two alternatives each have probability
$\tfrac12+o(1)$ along a subsequence, while the interval
length tends to infinity.

For $j\ge4$, choose integers $r_j\ge j^2$ increasingly.
With $r_j$ fixed, let $n_j$ be the least $n$ such that
\[
\Prob\bigl(\omega(G(n,1/j))\ge r_j\bigr)\ge\tfrac12.
\]
Choose $r_j$ sufficiently large that $n_j>n_{j-1}$ and
$a(n)\le1/j$ for every $n\ge n_j$.
Set $p(n)=1/j$ on $[n_j,n_{j+1})$, extending it constantly
to smaller $n$.

Define
\[
b_0=0,\qquad
b_r-b_{r-1}=3+\bigl|\{j\ge4:r_j\le r\}\bigr|.
\]
The interval lengths tend to infinity, so $b_r/r\to\infty$.
Within an interval, delete an isolated vertex;
at its first order, delete a clique vertex.
The deleted degree is at most $r-1=o(k)$ whenever
$b_{r-1}<k\le b_r$, verifying the deletion condition.
Moreover,
\[
b_{r_j}-b_{r_j-1}=j,
\qquad
b_{r_j}\le jr_j.
\]

Write $L_j=\log n_j$.
The first-moment bound gives
\[
\tfrac12\le n_j^{r_j}j^{-\binom{r_j}{2}},
\qquad
r_j=O(L_j/\log j)=o(L_j).
\]
Coupling by deletion of a uniformly chosen vertex shows
that the probability of containing an $r_j$-clique changes
by at most $r_j/n_j$ between $n_j-1$ and $n_j$.
Hence
\[
\tfrac12
\le\Prob\bigl(\omega(G(n_j,1/j))\ge r_j\bigr)
\le\tfrac12+\frac{r_j}{n_j}
=\tfrac12+o(1).
\]

Put $k_j=b_{r_j}\le jr_j$.
For a fixed $r_j$-clique, partition the remaining vertices
into $k_j-r_j$ blocks. Select one vertex from each block,
requiring no edges to the clique or previously selected vertices.
Each candidate succeeds with conditional probability at least
\[
(1-1/j)^{k_j}=n_j^{-o(1)}.
\]
A union bound over the blocks and all possible cliques
bounds the probability of any failure by
\[
n_j^{r_j}k_j
\exp\!\left\{-\frac{n_j}{2k_j}(1-1/j)^{k_j}\right\}
=o(1).
\]
Thus every $r_j$-clique whp lies in an admissible induced
subgraph of order $b_{r_j}$.

Consequently, for $G\sim G(n_j,1/j)$,
\[
\Prob(X_{\mathcal B}(G)\le b_{r_j-1})=\tfrac12+o(1),
\qquad
\Prob(X_{\mathcal B}(G)\ge b_{r_j})=\tfrac12+o(1).
\]
Since $b_{r_j}-b_{r_j-1}=j\to\infty$,
two-point concentration fails.
\end{proof}

\begin{example}[Sharpness of the condition $d(k)=o(k)$]
\label{ex:sharpness-d}
The class
\[
\mathcal A=\{H:|V(H)|\le2\alpha(H)\}
\]
satisfies the deletion condition with $d(k)=k-1$, but
$X_{\mathcal A}(G(n,1/2))$ is not two-point concentrated.
Thus $d(k)=o(k)$ cannot be replaced by $d(k)=O(k)$.
\end{example}

\begin{proof}
The construction doubles the independence number:
for every $n$-vertex graph $G$,
\[
X_{\mathcal A}(G)=\min\{n,2\alpha(G)\}.
\]
Indeed, writing $r=\alpha(G)$, every induced subgraph
$H\in\mathcal A$ satisfies
$|V(H)|\le2\alpha(H)\le2r$.
If $n\ge2r$, choose a set $S$ of $2r$ vertices containing a maximum independent set. Then
$\alpha(G[S])\ge r=|S|/2$, so $G[S]\in\mathcal A$.
If $n<2r$, then $G\in\mathcal A$.

For $G\sim G(n,1/2)$, we have $2\alpha(G)\le n$ whp,
and hence $X_{\mathcal A}(G)=2\alpha(G)$ whp.
Since a pair of consecutive integers contains only one
even integer, two-point concentration of $X_{\mathcal A}$ would imply one-point concentration of $\alpha(G(n,1/2))$,
contradicting the estimates of
Matula~\cite[Theorem~5]{Matula1976}.

Finally, the deletion condition holds with $d(k)=k-1$:
for even $k$, delete a vertex outside an independent set of size $k/2$; for odd $k$, delete any vertex.
Moreover, $K_{r,r}\in\mathcal A$ has all degrees equal to $r$, so no uniform sublinear degree bound is possible.
\end{proof}

\section*{Declaration of AI use}
The authors formulated the statement of Theorem~\ref{thm:main} as a conjecture and suggested using the deleted vertex's small degree in the overlap estimates.
The proof was then developed with the help of ChatGPT Pro 6; the authors take full responsibility for it.

\section*{Acknowledgements}

The research of Margarita Akhmejanova is supported by the Wenner-Gren
Foundations, grant number UPD2024-0154. The authors thank Maksim Zhukovskii for helpful comments.

\end{document}